\documentclass[12pt, epic, eepic]{amsart}
\usepackage{epsfig}

\newtheorem{theorem}{Theorem}[section]

\newtheorem{proposition}[theorem]{Proposition}

\newtheorem{remark}[theorem]{Remark}

\newtheorem{conjecture}[theorem]{Conjecture}

\newcommand{\Z}{\mathbb Z}

\newcommand{\Gal}{{\rm Gal}}

\newcommand{\Q}{\mathbb Q}
\newcommand{\comment}[1]{}

\title{Galois groups and an obstruction to principal graphs of subfactors}

\author{Marta Asaeda}

\address{Department of Mathematics, University of California Riverside,  900 Big Springs Drive, Riverside, CA, 92521 ,
 USA} 

\email{\tt marta@math.ucr.edu}

\thanks{The authors were sponsored in part by NSF grant
  \#DMS-0504199.}
  
\begin{document}
 
 \begin{abstract}
 The Galois group of the minimal polymonal of a Jones index
value gives a new type of obstruction to a principal graph, thanks
to a recent result of P.Etingof, D.Nikshych, and V.Ostrik.  We show that the sequence of the graphs given by Haagerup as candidates of principal graphs of subfactors, are not realized as principal graphs for $7<n \leq 27$ using GAP program. We further utilize Mathematica to extend the statement to $27 <n \leq 55$. We conjecture that none of the graphs are principal graphs for all $n>7$, and give an evidence using Mathematica for smaller graphs among them for $n>55$. The problem for the case $n=7$ remains open, however, it is highly likely that it would be realized as a principal graph, thanks to numerical computation by Ikeda. 
  \end{abstract}
\maketitle
 
 \section{Introduction} 
 
 Since V. F. R. Jones introduced the index theory of subfactors in \cite{J3}, the
theory of operator algebras have been achieving a remarkable development,
having relations with   low dimensional topology, solvable lattice model theory,
conformal field theory and quantum groups.  
 The following was the first breakthrough in the theory of Jones.
For a subfactor $N \subset M$, the index value $[M : N]$ belongs to the set
 $$ \{4 \cos^2 \frac{\pi}{n}  | \; n = 3, 4, 5 \cdots \} \cup [4, \infty].$$
 Later, he also introduced a principal graph and a dual principal
graph as finer invariants of subfactors. The Perron-Frobenius eigenvalue of the (dual) principal graph of a finite-depth subfactor $N \subset M$ is equal to $\sqrt{[M : N]}$. Using this fact 
  Jones proved in the middle of 1980's that subfactors with
index less than $4$ have one of the Dynkin diagrams of ADE type as their (dual)
principal graphs.     A. Ocneanu discovered a complete invariant for finite depth subfactors called ``paragroup" and announced
 in \cite{O2}
that subfactors with index less
than 4 are completely classified by the Dynkin diagrams $A_n$,
$D_{2n}$, $E_6$, and $E_8$. (See also \cite{BiN}, \cite{I1}, \cite{I3},
\cite{K2}, \cite{SV}.)  Ocneanu's invariant consists of numerical data called biunitary connection which is defined on squad of four graphs as depicted below:
\begin{figure}[h]
\begin{center}
\thinlines
\unitlength 1.0mm
\begin{picture}(30,20)(0,-7)
\multiput(11,-2)(0,10){2}{\line(1,0){8}}
\multiput(10,7)(10,0){2}{\line(0,-1){8}}
\multiput(10,-2)(10,0){2}{\circle*{1}}
\multiput(10,8)(10,0){2}{\circle*{1}}
\put(6,12){\makebox(0,0){$V_0$}}
\put(24,12){\makebox(0,0){$V_1$}}
\put(6,-6){\makebox(0,0){$V_2$}}
\put(24,-6){\makebox(0,0){$V_3$}}
\put(15,12){\makebox(0,0){$\mathcal G$}}
\put(6,3){\makebox(0,0){$\mathcal G$}}
\put(15,-6){\makebox(0,0){${\mathcal H}$}}
\put(24,3){\makebox(0,0){$\mathcal H$}}
\put(15,3){\makebox(0,0){$\alpha$}}
\end{picture} 
\end{center}
\caption{the biunitary connection $\alpha$ with four graphs}
\label{alpha}
\end{figure}
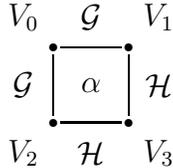
where Perron-Frobenius eigenvalues (PFEVs) of $\mathcal G$ and $\mathcal H$ coinside. When a biunitary connection satisfies an axiom called {\it flatness}, it is called a paragroup. Completeness of paragroups as invariants of subfactors is proved by S.Popa in \cite{P12}. In such case the graphs $\mathcal G$, $\mathcal H$ which are used to define a biunitary connection turns out to be the (dual) principal graph of the subfactor. Moreover each vertex of the graphs corresponds to a bimodule generated by the subfactor, and that the (dual) principal graph encodes the fusion rule of the bimodules.  Therefore, the problem of whether a given pair of graphs are realized as (dual) principal graph is reduced to construction of biunitary connection and to prove that it is flat. The detail of the paragroup theory is found in \cite{EK7}.  

\begin{figure}
$\psfig{figure=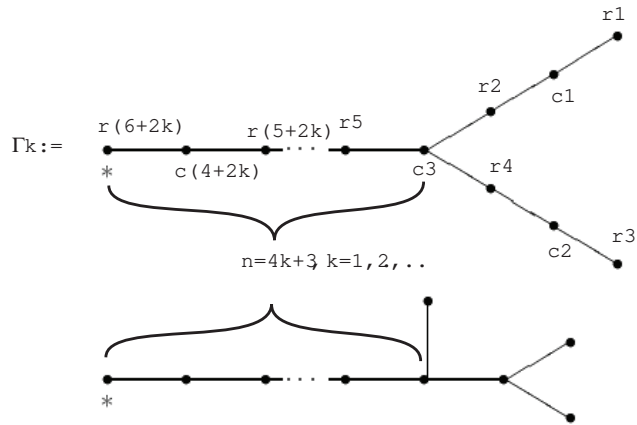,height=6cm}$
\caption{The the pairs of graphs (2) in the list of Haagerup}
\label{graph}
\end{figure}

After that, Popa (\cite{P12}) extended the correspondance between 
paragroups and subfactors of the hyperfinite II$_{1}$ factor to the 
strongly amenable case, and gave a
classification of subfactors with index equal to 4.  (In this
case, the dual principal graph of a subfactor is the same as the principal
graph.  See also \cite{IK}.) 
%%%%%%%%%%%
   Some subfactors with index larger than 4 had
been constructed from other mathematical objects.  For example, we can
construct a subfactor from an arbitrary finite group by a crossed product
with an outer action,
and this subfactor has an index equal to the order of the original
finite group.  Trivially, the index is at least 5 if it is larger than 4. 
We also have subfactors constructed from quantum groups
$U_q(sl(n)),\; q=e^{2 \pi i/k}$ with index $\frac{{\rm sin}^2 (n \pi 
/k)}{{\rm sin}^2(\pi/k)}$ as in \cite{Wn2} and these index values do not
fall in the interval $(4,5)$.  A subfactor with an index $3+\sqrt{3}=4.73...$ appeared in \cite{GHJ} is constructed by embedding the graph algebra of $A_{11}$ into that of $E_6$. 
U.Haagerup gave in 1991 a list of possible candidates of graphs
 which might be realized as
(dual) principal graphs of subfactors with index in 
$(4, 3+\sqrt{3})=(4, 4.732 \cdots)$ in \cite{H5}.  We see three pairs of finite graphs, including two pairs with parameters,
along with a pair of infinite graphs $A_\infty$ in (1),  in \S7 of \cite{H5}.  
  Since then  D. Bisch  proved that a subfactor with (dual) principal graph (4) in \S7 of \cite{H5} does {\it not} exist \cite{Bs8} by checking inconsistency of fusion rule on the graph. Haagerup and the author proved that two pairs of graphs:   the case 
 $n=3$ of (2): see Figure \ref{graph}
  as well as the case (3) in \S7 of \cite{H5}, are realized as (dual) principal graphs of subfactors, and that such subfactors are unique respectively (\cite{AH}). 
 The remaining problem was whether the graphs for the case $n >3$ of (2) as in  Figure \ref{graph} would be realized as (dual) principal graphs of subfactors. Haagerup proved that the obstruction as found for the case (4) by Bisch does not exist on any of the pairs of the graph in (2). Moreover he proved that unique biunitary connection exists for each the pair of the graphs (\cite{Hp}). 
 For the case $n=7$, it was numerically checked by K.Ikeda that the biunitary connection is flat(\cite{Ik}). 
 There has been no progress for several years since then. Recently P.Etingof, D.Nikshych, and V.Ostrik showed in \cite{ENO} Theorem 8.51, that the index of a subfactor has to be a cyclotomic integer, namely an algebraic integer that lies in a cyclotomic field. This implies that if the square of PFEV of a graph is not a cyclotomic integer, the graph cannot be the (dual) principal graph of a subfactor. 
 
 In this paper we prove that the graphs in  Figure \ref{graph} are not (dual) principal graph for $n=4k+3$ for $1<k \leq 13$ by computing Galois groups of minimal polynomials of the square of PFEV of the graphs. We also prove that for the case $k=1$ the square of PFEV of the graphs is a cyclotomic integer: thus we cannot eliminate the possibility that the graphs might be (dual) principal graphs of a subfactor.  We further give an evidence by Mathematica computation that the graphs for larger $k$ are not principal graphs. We conjecture that for the case $k >1$ none of the graph is (dual) principal graph of a subfactor. 
     
The author is much grateful to Y. Kawahigashi, D.Bisch, and V.Jones  for informing me of the result in \cite{ENO} and suggesting the possibility of utilizing Galois theory. She also thanks U.Haagerup for sharing with me his unpublished results on the larger graphs, and L.Washington, T.Saito for  helpful discussion on Galois theory. 
% \subsection{Galois theory}
\section{Preliminaries}
In the following we list some known theorems in Galois theory necessary for later discussion.   
\begin{theorem}(Galois correspondence  \cite{Hun}) \\
Let $F$ be a finite dimensional Galois extension of $K$. Then there is one-to-one correspondence between the set of all intermediate fields of the extension and the set of all subgroups of the Galois group $G={\rm Aut}_K F$, given by $E \mapsto {\rm Aut}_E F$.  An intermediate field $E$ is Galois over $K$ if and only if $ {\rm Aut}_E F$ is a normal subgroup of $G$. In this case $Aut_K E \cong G/{{\rm Aut}_E F}$. 
\end{theorem}
\begin{theorem} (Kronecker-Weber, \cite{Wash}) \\
\label{kron}
Let $E$ be an extension of $\Q$ so that Aut$_\Q E$ is abelian. Then the field $E$ is a subfield of a cyclotomic field, that is, a field extension of $\Q$ by a primitive root of unity. 
\end{theorem}
\begin{theorem} (\cite{Mil})
\label{discriminant}
Let $p(x) \in \Q[x]$ be a monic irreducible polynomial of degree $n$, and let $G$ be its Galois group. Let $\xi_1, ... \xi_n$ be the roots of $p$. We define the discriminant of $p$ by $D_p={\Delta_p}^2$, where
$$\Delta_p= \prod_{i < j} (\xi_i-\xi_j)$$
up to sign.
Then we have
$${\Delta_p}^2 \in \Q,$$
and 
$$ G \subset \frak{A}_n \Longleftrightarrow { \Delta_p} \in \Q.$$
We may replace $\Q$ by $\Z$ if $p(x) \in \Z[x]$.
\end{theorem}
\begin{remark} (\cite{vdW}, \S 26-28)
\label{discriminantformula}
The discriminant $D_p$ of a monic polynomial $p(x)$ relates to the resultant Res$(p,p')$ of $p$ and $p'$  by $D_p=(-1)^{n(n-1)} {\rm Res}(p,p')$.    (Note that in \cite{vdW} the sign is omitted.) The resultant Res$(p,p')$ is given as follows: Let $p(x)=x^n + a_{n-1} x^{n-1} + ... a_1 x + a_0$. Then the resultant Res$(p,p')$ is equal to
$$
 \left|\begin{array}{cccccccccc}1 & a_{n-1} & a_{n-2} & \cdots & \cdots & a_0 & 0 & \cdots & 0 & 0 \\0 & 1 & a_{n-1} & a_{n-2} & \cdots & \cdots & a_0 & 0 & \cdots & 0 \\
 \vdots & \ddots & \ddots & \ddots & \ddots & \ddots & \ddots & \ddots & \ddots & \vdots \\
  \vdots & \ddots & \ddots & \ddots & \ddots & \ddots & \ddots & \ddots & \ddots & 0 \\
 0 & \cdots & \cdots & 0  & 1 & a_{n-1} & a_{n-2} & \cdots & \cdots & a_0  \\
 n & b_{n-1} & b_{n-2} & \cdots & \cdots & b_0 & 0 & \cdots & 0 & 0 \\0 & n & b_{n-1} & b_{n-2} & \cdots & \cdots & b_0 & 0 & \cdots & 0 \\
 \vdots & \ddots & \ddots & \ddots & \ddots & \ddots & \ddots & \ddots & \ddots & \vdots \\
  \vdots & \ddots & \ddots & \ddots & \ddots & \ddots & \ddots & \ddots & \ddots & 0 \\
 0 & \cdots & \cdots &  0 & n & b_{n-1} & b_{n-2} & \cdots & \cdots & b_0
 \end{array}\right|,
 $$
 where $b_i:=ia_i$. 
 It is easily obtained that  $D_p=(-1)^c |D_p|$, where $c$ is half the number of complex roots of $p(x)$. 
 \end{remark}
 
 Our use of these theorems are as follows: Let $d$ be an algebraic integer with the minimal polynomial $p(x) \in \Z[x]$. Suppose $d$ lies in a cyclotomic field $F$. Then $\Q(d)$ is an intermediate field of $F/\Q$. Thus by Galois correspondence it corresponds to a subgroup of Aut$_\Q F=\Z_n$ for some $n$. Since any subgroup of $\Z_n$ is normal, the extension $\Q(d)/\Q$ is Galois. Thus all the roots of $p(x)$ belong to $\Q(d)$, and the Galois group of $p(x)$ is a quotient of $\Z_n$, which is always abelian. Together with Kronecker's theorem, we conclude the following:
 \begin{proposition}
 An algebraic number $d$ is a cyclotomic integer if and only if its minimal polynomial has an abelian Galois group. 
 \end{proposition}
 
 The following facts are useful in computing Galois group. 
 \begin{proposition} (\cite{Mil})
 \label{transitive}
Let $K$ be a field,  $f(x) \in K[x]$, and $E$ be a splitting field of $f$ generated by its roots over $K$. If $f$ is irreducible and separable, $\Gal(E/K)=\Gal(f)$ acts transitively on the roots of $f$, i.e.   $\Gal(f)$  is a transitive subgroup of $\frak{S}_n$, where $n={\rm deg}f$. Furthermore $|\Gal(f)|$ is divisible by $n$. 
%If $f$ is not irreducible, roots that belong to distinct irreducible component are not maped to each other by an element of $\Gal(E/K)$. 
 \end{proposition}
 %
% \begin{proposition}
% Let $K$ be a finite field. Then Galois group of any irreducible polynomial in $K[x]$ is cyclic. 
% \end{proposition}
% %
% \begin{proposition} ( Dedekind \cite{vdW})
% \label{cycle}
% Let $f(x) \in \Z[x]$ with deg$f=n$ be an irreducible monic polynomial. For a prime number $p$, let $\bar{f}(x) \in \Z_p[x]$ be the image of $f$ by the natural homomorphism induced by $\Z \to \Z_p$, and let 
% $$\bar{f}(x)=\bar{f}_1(x)\cdots \bar{f}_i(x)$$
% be the irreducible factorization with deg$\bar{f}_j=n_j$.  Suppose that all $\bar{f}_j$'s are distinct in $\Z_p[x]$. Then $\Gal_\Q(f)$ contains a permutation of the type $(1 \ 2 \ \cdots \ n_1)(n_1+1 \ \cdots \ n_1+n_2)\cdots=:[n_1,n_2, ... , n_i]$. 
% \end{proposition}
% 

 %
 %
 %
\section{Minimal polynomials for the square of Perron-Frobenius eigenvalues}
 In this section we give a formula for the polynomials which are candidates for minimal polynomials of the square of PFEVs of the graphs in Figure \ref{graph}. Since PDEVs of each pair coincide, we only use the first sequence of the graphs. The adjacency matrix of the graph $\Gamma_k$ is as follows:
 $$A_k:=
\bordermatrix{ &  c_1 &c_2   &c_3   &c_4   &\cdots   &  \cdots &  c_{4+2k}   \cr
r_1&1 & 0 & 0 & 0 & \cdots & \cdots & 0    \cr
r_2&1 & 0 & 1 &0& \cdots  & \cdots  & 0&   \cr
r_3& 0 & 1 & 0 & 0& \cdots &   \cdots & 0  \cr
r_4& 0 & 1 & 1 & 0 & \cdots & \cdots & 0    \cr
r_5& \vdots & 0 & 1 & 1 & 0 &  \cdots & 0    \cr
\vdots & \vdots &   \vdots &   \ddots & \ddots & \ddots & \ddots & \vdots    \cr
r_{5+2k}  & 0 & 0& \cdots & 0 & 1 & 1 & 0    \cr
r_{6+2k} & 0 & 0& \cdots & \cdots & 0 & 1 & 1    
  },
  $$
Where $(i, j)$-entry is given by the number of edges connecting $r_i$ and $c_j$. Notice that $r_i$'s are even vertices and $c_j$'s are odd vertices of the graph $\Gamma_k$, considering $*$ as zero-th vertex.   The PFEV $\beta_k$ is the unique eivenvalue of the matrix 
  $$M_k:=\left[ \begin{array}{cc} 
  {\bf 0}_{6+2k} & A_k \\
  {A_k}^t & {\bf 0}_{4+2k}
  \end{array} \right]
  $$
  with the largest norm. Note that the rows (resp. columns) of $M_k$ are labeled by the vertices of $\Gamma_k$ in the order of $r_1, r_2, ... , r_{6+2k}, c_1, ... , c_{4+2k}.$ It is known that $\beta_k$ is a real number, its eigenspace is one dimensional, and an eigenvector can be taken to be a real vector. Let $u_k$ be the eigenvector chosen so that the entry corresponding to the row labeled by $r_{6+2k}$ will be one. We may regard $u_k$ as a direct sum of two vectors $u_k=(v_k, w_k)$, where $v_k$ consists of the entries corresponding to even vertices (i.e. $r$'s), $w_k$ corresponds to odd vertices (i.e. $c$'s). Then we have the following relation: 
  \begin{eqnarray*}
  A_k w_k &=&\beta_k v_k, \\
  (A_k)^t v_k &=&  \beta_k w_k. 
  \end{eqnarray*}
   Consider
 $$ (M_k)^2= \left[ \begin{array}{cc} 
   A_k {A_k}^t& {\bf 0}  \\
{\bf 0}  & {A_k}^t A_k
  \end{array} \right]. 
  $$
  The largest norm of the eigenvalues is given by ${\beta_k}^2=:d_k$. On the other hand clearly $d_k$ is an eigenvalue with eivenvector $u_k$, thus $d_k$ is the PFEV of 
  $(M_k)^2$. Since non-zero eigenvalues of $A_k {A_k}^t$ and that of ${A_k}^t A_k$ coincides, $d_k$ must be the PFEV of  ${A_k}^t A_k$ (resp.  $A_k {A_k}^t$). Thus we deal with $N_k:={A_k}^t A_k$, since it is a smaller matrix. $N_k$ is given as follows:
$$
  N_k=
\bordermatrix{ &  c_1 &c_2   &c_3   &c_4      &\cdots&\cdots&\cdots &   c_{3+2k} &  c_{4+2k}   \cr
 c_1 &2 & 0 & 1 & 0 & \cdots & \cdots & \cdots & \cdots & 0       \cr
 c_2 &0 & 2 & 1 & 0 &   &   &   &   & \vdots       \cr
 c_3 & 1 & 1 & 3 & 1 & 0 &   &   &   & \vdots     \cr
 c_4  & 0 & 0 & 1 & 2 & 1 & 0 &   &   & \vdots       \cr
c_5& 0 & 0 & 0 & 1 & 2 & 1 & 0 &   & \vdots      \cr
 \vdots & \vdots &   &   & \ddots & \ddots & \ddots & \ddots & \ddots & \vdots      \cr
c_{4+2(k-1)} & \vdots &   &   &   & 0 & 1 & 2 & 1 & 0      \cr
 c_{3+2k} & 0 & \cdots & \cdots & \cdots & \cdots & 0 & 1 & 2 & 1     \cr
 c_{4+2k} & 0 & \cdots & \cdots & \cdots & \cdots & \cdots & 0 &  1 & 2 } 
$$
Let us call the $(3+2k) \times (3+2k)$ submatrix of $N_k$ containing first $3+2k$ columns and rows $N_{k-1/2}$. Let $p_k$ be the characteristic polynomial of $N_k$. Then clearly for every half integer $k > 1$ we have the following recursive relation:
$$ p_k(x)= (2-x) p_{k-1/2}(x) - p_{k-1}(x),$$ 
and from this we easily get 
\begin{eqnarray*}
p_k(x)&=&(x^2-4x+2)p_{k-1}(x)-p_{k-2}, \\
p_0(x)&=&(x^2-5x+3)(x-2)^2, \\
 p_1(x)&=&(x^3-8x^2+17x-5)(x-2)^2(x-1).
\end{eqnarray*}
Since both $p_0$ and $p_1$ contain $(x-2)^2$ as factor, $p_k(x)$ for any $k$ also contains $(x-2)^2$. Since our concern is minimal polynomials for PFEVs, we set $q_k(x):=p_k(x)/(x-2)^2$. Obviously $q_k(x)$ satisfy the same recursive equation. Furthermore we have the following:
\begin{proposition}
The polynomial $q_k(x)$ is divisible by $(x-1)$ if and only if $k \equiv 1$ mod $3$. 
\end{proposition}
\noindent
{\bf Proof}\\
We simply plug in $x=1$  in the equation: we know that $q_0(1)=-1$, $q_1(1)=0$, and we easily obtain $q_2(1)=1$.  Suppose $q_{k-2}(1)=0$   Then $q_k(1)=-q_{k-1}(1)$, thus $q_{k+1}(1)=-q_k(1)-q_{k-1}(1)=0$, i.e. $q_{k+1}(x)$ is divisible by $(x-1)$ in this case. Suppose  $q_{k-2}(1)=\pm 1$. Then $ q_k(1)= -q_{k-1}(1) \mp 1$, thus $q_{k+1}(1)=-q_k(1)-q_{k-1}(1)=\pm 1.$ \qed 

\bigskip
We may solve the recursive equation explicitly in a standard method: we obtain 
  the following:
$$q_k(x)=A(x) a(x)^{2k} + B(x) b(x)^{2k},$$
where $a(x) = (2 - x + \sqrt{x^2 - 4x})/2$, $b(x) = (2 - x -\sqrt{x^2 - 4x})/2$, $A(x)=\frac{-1}{a(x)^2-b(x)^2}(q_0(x)b(x)^2-q_1(x))$, and $B(x)=\frac{1}{a(x)^2-b(x)^2}(q_0(x)a(x)^2-q_1(x))$. We conjecture the following:
\begin{conjecture}
Let 
$$
r_k(x)= \left\{ \begin{array}{c} 
q_k(x)/(x-1), \; {\rm if} \; k \equiv 1 \ {\rm mod}\  3, \\
q_k(x), \; {\rm else}.  
\end{array}
\right.
$$
Then $r_k(x)$ is irreducible for any $k$. 
\end{conjecture}
This is  proved for the values of $k$ up to $6$ by using GAP \cite{Gap}. (For $k=0,1$ it is obvious.) GAP is an open-source program designed for mathematicians that does not give any approximate solutions without an message stating so or deliberately set by each user to give an approximation, i.e. any solutions given by GAP is mathematically accurate at least to the level of published results.  For larger $k$, it is checked by Mathematica for individual case.   Mathematica is not an open-source software, 
and the users have no way of knowing the reliability of results, thus we do not dare to claim it as a proof. However we would like to note that we found a prime number for each $k$ modulo which Mathematica thinks that $r_k(x)$ is irreducible, for $k=7, ... , 13.$ Following is the list of the smallest primes $p$ used for each $k$: we list as $(k,p)$. \\
(7,3), (8,2), (9,5), (10,3), (11,3), (12, 2), (13, 11).  \\
Since factorization of a polynomial modulo prime is a finite process, we can conclude that $r_k(x)$'s for $k=7, ..., 13$ are irreducible. And for $k=14, ... , 19$,  Mathematica thinks  that $r_k(x)$'s are irreducible for whatever reason we do not know. 
 Note that the number ``19" is totally arbitrary: it does not mean that it is the maximum $k$ that Mathematica could handle. 
\section{Galois groups of $r_k(x)$ and cyclotomicity of $d_k$.}
Our aim is to check if $d_k$ is cyclotomic number. Thus we need to compute the Galois group of its minimal polynomial $m_k(x)$. If $r_k(x)$ in the previous section is irreducible, it coincides with $m_k(x)$. Otherwise $m_k(x)$ factors $r_k(x)$. 
For $k=0,1$, we have the following:
\begin{proposition}
\begin{eqnarray*}
\Gal(r_0(x))&=&\Z_2 \\
\Gal(r_1(x))&=&\Z_3, 
\end{eqnarray*}
where for a polynomial $f \in \Q[x] $, we denote its Galois group by $Gal(f)$. 
\end{proposition}
\noindent
{\it Proof} \\
Using Proposition \ref{transitive}, we have $\Gal (r_0(x)) = \Z_2=\frak{S}_2$ and that $\Gal (r_1(x))$ is a transitive subgroup of  $\frak{S}_3$, i.e. $\frak{A}_3=\Z_3$ or $\frak{S}_3$. Now, we compute the discriminant of $r_1(x)=x^3-8x^2+17x-5$ using the formula given in Remark \ref{discriminantformula}.  We have the resultant Res$(r_1,r'_1)=-139$, thus  
$D_{r_1}= (-1)^{3 \cdot 2/2} \cdot  (-169)=13^2$. Therefore $\Gal(r_1(x))=\Z_3. $ \qed

The above result implies that $d_0$, $d_1$ are cyclotomic numbers by Theorem \ref{kron}. For $k=0$ in fact we had already known that $d_0$ is cyclotomic, since the graph $\Gamma_0$ in Figure \ref{graph} is realized a principal graph. For $k=1$, this result implies that $\Gamma_1$  still has a chance of being a principal graph, but it needs to be checked by other methods. 

 We may still utilize GAP for small values of $k$.  
 \begin{proposition}
 For $k=2,..., 6$, $\Gal(r_k)=\frak{S}_{n_k}$, where $n_k={\rm deg}(r_k)$. 
 \end{proposition}
 This implies that $d_k$'s for $k=2,..., 6$ are not cyclotomic, thus corresponding $\Gamma_k$'s cannot be realized as principal graphs of subfactors. The readers may wish to check it on their own for small $k$, here we provide the first two polynomials: $r_2(x)=x^6 - 13 x^5 + 63 x^4 - 140 x^3 + 142x^2 - 59x + 7$, $r_3(x)=x^8- 17x^7 + 117 x^6-418 x^5+ 827x^4- 898x^3+ 502x^2 - 124x +9$. 

 For larger $k$, we may still work using Mathematica. We have the following very strong fact: \begin{proposition} (\cite{Kon})
 Let $f(x) \in \Z[x]$ be an irreducible polynomial with degree $n$. Then $\Gal (f)= \frak{S}_n$ if the discriminant of $f$ is square-free. 
 \end{proposition}
 This in particular implies the following: 
 \begin{theorem}
 Let $\Gamma$ be a finite graph, and $d_{\Gamma}$ be the square of PFEV of $\Gamma$, and let $m_\Gamma(x)$ its minimal polynomial with degree larger than $2$.  Then $\Gamma$ is not realized as a principal graph of a subfactor if the discriminant of $m_\Gamma$ is square-free. 
 \end{theorem}
 Mathematica does not have a command for discriminant, however it does have resultant, thus we may compute discriminant up to sign.  Below we show the mathematica computations. $ired[k,x]$ corresponds to $r_{k-1}(x)$. D[{\it function}, x] is the derivative of a given function in $x$.  
 \begin{eqnarray*}
\rm{dd[k\_]} &:=& {\rm{Abs[Resultant[Simplify[ired[k, x]], D[Simplify[ired[k, x]], x], x]]}}  \\
&& {\mbox{(discriminant of }} r_{k+1}(x))  {\mbox{up to sign) }}\\
{\rm{fd[k\_]}} &:=& {\rm{FactorInteger[Abs[Resultant[Simplify[
    ired[k, x]], }} \\ && {\rm{D[Simplify[ired[k, x]], x], x]]] }}
  \ \ {\mbox{ (find  the factorization of } }{\rm dd[k])} \\
\end{eqnarray*} 
The followings are the results: note again that the numbering is shifted by $1$ from our paper. Since the computation gives integers as result, we can take the computation to be accurate. 
Here  $\{\{p,n_p\},\{q,n_q\}, ...\} $ means that the given number is prime-factorized in the form of $p^{n_p} \cdot q^{n_q} \cdot ....$. 
\begin{eqnarray*}
&& {\rm fd[3]} = \\ 
   &&\{\{1471, 1\}, \{5171, 1\}\} \\
&& {\rm fd[4]} = \\ 
   &&\{\{1097, 1\}, \{4261, 1\}, \{8677, 1\}\} \\
&& {\rm fd[5]} = \\ 
   &&\{\{281, 1\}, \{643, 1\}, \{31281527, 1\}\} \\
&& {\rm fd[6]} = \\ 
   &&\{\{192667, 1\}, \{47117433796403, 1\}\}\\
 && {\rm fd[7]} = \\ 
   && \{\{3, 1\}, \{47, 1\}, \{3323, 1\}, \{3613, 1\}, \{7487, 1\}, \{22182696017, 1\}\}  \\
 &&{\rm fd[8]} 
 = \\  && \{\{29, 1\}, \{1427, 1\}, \{11933, 1\}, \{35419, 1\}, \{595801, 1\}, \{7143737, 1\}\} \\
 &&{\rm fd[9]} 
 = \\ &&  \{\{769765583537031753607466863873613, 1\}\} \\
  &&{\rm fd[10]} 
 = \\  &&  \{\{31, 1\}, \{1625255809, 1\}, \{1226665686533457543318366623, 1\}\} \\
  &&{\rm fd[11]} 
 = \\ &&  \{\{230113699, 1\}, \{1990348035579493, 1\}, \{47658861361724611, 1\}\} \\
   &&{\rm fd[12]} 
 = \\ &&   \{\{119813, 1\}, \{23296847792041232351, 1\}, \{285981481927230196531187, 1\}\} \\
 &&{\rm fd[13]} 
 = \\ &&   \{\{17761, 1\}, \{15894547, 1\}, \{202995484303, 1\},  \\
 && \{    993013213822241, 1\}, \{2155998286155473, 1\}\} \\
 &&{\rm fd[14]} 
 = \\  &&  \{\{745621, 1\}, \{21802562773909, 1\}, \{468985859471443, 1\}, \\
&&\{6697788892778259550891, 1\}\} \\
&&{\rm fd[15]} 
 = \\ &&   \{\{41, 1\}, \{85503853,
   1\}, \{525628115273, 1\},  \\
   && \{2640893817458692409629597355203029150324817, 1\}\} 
\end{eqnarray*}
\begin{eqnarray*}
   &&{\rm fd[16]} 
 = \\  &&   \{\{3, 1\}, \{943618253, 1\}, \{9374569646597215017911, 
  1\},  \\
   &&  \{46088050874425115317503501160573626593, 1\}\} \\
    &&{\rm fd[17]} 
 = \\ &&   \{\{1487, 1\}, \{
    14737, 1\}, \{42895179574588531, 1\}, \\
&& \{602237386867482390429552519214023674351054569169, 1\}\} \\
  &&{\rm fd[18]} 
 = \\  &&  \{\{281, 1\}, \{619, 1\}, \{21149, 1\}, \{2454047,
   1\}, \{27050115645481, 1\},  \\ && \{109185979881289, 
    1\}, \{4378070972266731488149874970904128697, 1\}\} \\
&&{\rm fd[19]} 
 = \\ &&   \{\{408019, 1\}, \{1085473, 1\}, \{31856719917482639623, 1\}, \\ 
 && \{  3088206373486625469392445666425929334732872642867424521, 1\}\} \\
 &&{\rm fd[20]} 
 = \\ &&  \{\{29, 1\}, \{4783, 1\}, \{4700047160321, 1\}, \\ 
 && \{
  3332315546190627198162521685721451274758792227246905 \\
  && 7938406361694282211, 1\}\}
\end{eqnarray*}
Observe that all the exponent of the prime factors are $1$, namely they are all square-free. Thus we conclude the following:
\begin{theorem}
The graphs $\Gamma_k$ in Figure \ref{graph} are not principal graphs of subfactor for  $13  \geq k > 1$. 
\end{theorem}
Note that from discriminant test it is also very likely that the theorem is true for the case $19 \geq k>13$: the only reservation for this range is irreducibility of $r_k$'s, as noted earlier. 
 At this point we hope for a proof for the following conjecture: 
\begin{conjecture}
The graphs $\Gamma_k$ in Figure \ref{graph} are not principal graphs of subfactor for any $k > 1$. 
\end{conjecture}

\end{document}